\documentclass[a4paper]{amsart}

\usepackage{amsmath,amssymb,amsthm}

\usepackage{hyperref}
\hypersetup{colorlinks=true,hyperindex,plainpages=false,
            pdftitle={Gauge Symmetries and Noether Currents in Optimal Control},
            pdfkeywords={optimal control, Pontryagin extremals, gauge symmetry,
                         second Noether theorem, Noether currents.},
            pdfauthor={Delfim F. M. Torres},
            pdfcreator={LaTeX with package hyperref}
}

\newtheorem{theorem}{Theorem}[section]

\theoremstyle{remark}
  \newtheorem{remark}{Remark}[section]

\theoremstyle{definition}
  \newtheorem{definition}{Definition}[section]

\renewcommand\MR[1]{\href{http://www.ams.org/mathscinet-getitem?mr=#1}{\textsf{MR}~#1}}
\newcommand\ZBL[1]{\href{http://www.emis.de/MATH-item?#1}{\textsf{Zbl}~#1}}
\newcommand\JFM[1]{\href{http://www.emis.de/cgi-bin/JFM-item?#1}{\textsf{JFM}~#1}}

\begin{document}

\title[Gauge Symmetries in Optimal Control]{Gauge Symmetries
and Noether Currents in Optimal Control}

\author{Delfim F. M. Torres}

\address{Departamento de Matem\'{a}tica,
Universidade de Aveiro, 3810-193 Aveiro, Portugal}

\email{\href{mailto:delfim@mat.ua.pt}{delfim@mat.ua.pt}}

\urladdr{\href{http://www.mat.ua.pt/delfim}{http://www.mat.ua.pt/delfim}}

\thanks{This research was supported in part by the
\emph{Optimization and Control Theory Group} of
the R\&D Unit \emph{Mathematics and Applications},
and the program PRODEP III 5.3/C/200.009/2000.
Partially presented at the 5th Portuguese Conference on
Automatic Control (\emph{Controlo 2002}), Aveiro, Portugal, September 5--7, 2002.
Accepted for publication in \emph{Applied Mathematics E-Notes}, Volume 3.}

\date{}

\subjclass[2000]{49K15, 49S05.}

\keywords{optimal control, Pontryagin extremals, gauge symmetry,
second Noether theorem, Noether currents.}


\begin{abstract}
We extend the second Noether theorem to optimal control problems
which are invariant under symmetries depending upon
$k$ arbitrary functions of the independent variable and their derivatives
up to some order $m$. As far as we consider a semi-invariance notion,
and the transformation group may also depend on the control variables,
the result is new even in the classical context of the calculus of variations.
\end{abstract}


\maketitle


\section{Introduction}

The study of invariant variational problems
\begin{equation*}
\text{Minimize } J\left[x(\cdot)\right] = \int_a^b L\left(t,x(t),\dot{x}(t)\right)
\mathrm{d}t
\end{equation*}
in the calculus of variations
was initiated in the early part of the XX century by Emmy Noether who,
influenced by the works of Klein and Lie on the transformation properties
of differential equations under continuous groups of transformations
(\emph{see} \textrm{e.g.} \cite[Ch.~2]{MR2001f:35024}), published
in her seminal paper \cite{JFM46.0770.01,MR53:10538} of 1918
two fundamental theorems, now classical results and known as the (first) Noether
theorem and the second Noether theorem, showing that invariance with respect
to a group of transformations of the variables $t$ and $x$
implies the existence of certain conserved quantities.
These results, also known as \emph{Noether's symmetry theorems},
have profound implications in all physical theories,
explaining the correspondence between symmetries of the
systems (between the group of transformations acting
on the independent and dependent variables of the system)
and the existence of conservation laws.
This remarkable interaction between the concept of invariance
in the calculus of variations and the existence of first integrals
(Noether currents) was clearly recognized by Hilbert \cite{hilbert}
(\textrm{cf.} \cite{MR57:13703}).

The first Noether theorem establishes the existence of $\rho$
first integrals of the Euler-Lagrange differential equations when
the Lagrangian $L$ is invariant under a group of
transformations containing $\rho$ parameters.
This means that the invariance hypothesis leads to quantities which are
constant along the Euler-Lagrange extremals.
Extensions for the Pontryagin extremals of optimal
control problems are available in
\cite{delfim3ncnw,delfimEJC,torresCM02I01E}.

The second Noether theorem establishes the existence of $k \left(m + 1\right)$
first integrals when the Lagrangian is invariant under
an infinite continuous group of transformations
which, rather than dependence on
parameters, as in the first theorem,
depend upon $k$ arbitrary functions and their derivatives up to
order $m$. This second theorem is not as well known as the first.
It has, however, some rather interesting implications.
If for example one considers the functional of the basic problem of the calculus
of variations in the autonomous case,
\begin{equation}
\label{e:PBCV:autonomo}
J\left[x(\cdot)\right] = \int_a^b L\left(x(t),\dot{x}(t)\right) \mathrm{d}t \, ,
\end{equation}
the classical Weierstrass necessary optimality
condition can easily be deduced from the fact that the integral
\eqref{e:PBCV:autonomo} is invariant under transformations of the form
$T = t + p(t)$, $X = x(t)$, for an arbitrary function $p(\cdot)$
(see \cite[p.~161]{MR58:18024}).
The second Noether theorem is related to: (i) parameter invariant
variational problems, \textrm{i.e.},
problems of the calculus of variations, as in the homogeneous-parametric form,
which are invariant under arbitrary transformations of the independent variable $t$
(see \cite[p.~266]{MR16:426a}, \cite[Ch. 8]{MR58:18024},
\cite[p.~179]{ZBL0964.49001}); (ii) the singular
Lagrangians and the constraints in the Hamiltonian formalism,
a framework studied by Dirac-Bergmann (see \cite{MR2001k:70023,guerraPhD});
(iii) the physics of gauge theories,
such as the gauge transformations of electrodynamics,
electromagnetic field, hydromechanics, and relativity
(see \cite[pp. 186--189]{ZBL0964.49001},
\cite[p. 160]{MR58:18024}, \cite{MR50:14419}, \cite{MR53:10537}).
For example, if the Lagrangian $L$ represents
a charged particle interacting with a electromagnetic field, one finds
that it is invariant under the combined action of the so called gauge
transformation of the first kind on the charged particle field,
and a gauge transformation of the second kind on the electromagnetic field.
As a result of this invariance it follows, from second
Noether's theorem, the very important conservation of charge.
The invariance under gauge transformations is a basic requirement
in Yang-Mills field theory, an important subject, with many questions
for mathematical understanding (\textrm{cf.} \cite{CMIYangMill}).

To our knowledge, no second Noether type theorem is available for the
optimal control setting. One such generalization is our concern here.
Instead of using the original argument \cite{JFM46.0770.01,MR53:10538}
of Emmy Noether, which is fairly complicated
and depends on some deep and conceptually difficult results
in the calculus of variations, our approach follows,
\emph{mutatis mutandis}, the paper \cite{delfim3ncnw}, where
the first Noether theorem is derived almost effortlessly by means of elementary
techniques, with a simple and direct approach, and it is motivated by the
novelties introduced by the author in \cite{delfimEJC}.
Even in the classical context (\textrm{cf. e.g.} \cite{MR50:14419})
and in the simplest possible situation,
for the basic problem of the calculus of variations,
our result is new since we consider symmetries of the system which alter
the cost functional up to an exact differential; we introduce
a semi-invariant notion with some weights $\lambda^0,\ldots,\lambda^m$
(possible different from zero);
and our transformation group may depend also on $\dot{x}$ (the control).
Our result hold both in the normal and abnormal cases.


\section{The Optimal Control Problem}

We consider the optimal control problem in Lagrange form
on the compact interval $[a,b]$:
\begin{equation*}
\text{Minimize } J[x(\cdot),u(\cdot)] =
\int_a^b L\left(t,x(t),u(t)\right) \mathrm{d}t
\end{equation*}
over all admissible pairs $\left(x(\cdot),u(\cdot)\right)$,\footnote{The notation
$W_{1,1}$ is used for the class of absolutely continuous functions, while
$L_{\infty}$ represents the class of measurable and
essentially bounded functions.}
\begin{equation*}
\left(x(\cdot),u(\cdot)\right) \in
W_{1,1}^n\left([a,b];\mathbb{R}^n\right) \times
L_{\infty}^r\left([a,b];\Omega \subseteq \mathbb{R}^r\right) \, ,
\end{equation*}
satisfying the control equation
\begin{equation*}
\dot{x}(t) = \varphi\left(t,x(t),u(t)\right) \quad \text{a.e. } t \in [a,b] \, .
\end{equation*}
The functions
$L : \mathbb{R} \times \mathbb{R}^n \times \mathbb{R}^r \rightarrow \mathbb{R}$
and
$\varphi : \mathbb{R} \times \mathbb{R}^n \times \mathbb{R}^r \rightarrow \mathbb{R}^n$
are assumed to be $C^1$ with respect to all variables
and the set $\Omega$ of admissible values of the control parameters is
an arbitrary open set of $\mathbb{R}^r$.

Associated to the optimal control problem there is the Pontryagin Hamiltonian
$H : [a,b] \times \mathbb{R}^n \times \Omega \times \mathbb{R}
\times \left(\mathbb{R}^n\right)^T \rightarrow \mathbb{R}$
which is defined as
\begin{equation}
\label{eq:Hamiltonian}
H(t,x,u,\psi_0,\psi) = \psi_0 L(t,x,u) + \psi \cdot \varphi(t,x,u) \, .
\end{equation}
A quadruple $\left(x(\cdot),\,u(\cdot),\,\psi_0,\,\psi(\cdot)\right)$,
with admissible $\left(x(\cdot),\,u(\cdot)\right)$,
$\psi_0 \in \mathbb{R}^{-}_{0}$, and $\psi(\cdot) \in
W_{1,1}\left([a,b];\,\mathbb{R}^n \right)$
($\psi(t)$ is a covector $1 \times n$), is called
a \emph{Pontryagin extremal}  if the following two
conditions are satisfied for almost all $t \in [a,b]$:
\begin{description}
\item[The Adjoint System]
\begin{equation}
\label{eq:sh}
\dot{\psi}(t) = - \frac{\partial H}{\partial
x}\left(t,\,x(t),\,u(t),\,\psi_0,\,\psi(t)\right) \, \text{;}
\end{equation}
\item[The Maximality Condition]
\begin{equation}
\label{eq:cm}
H\left(t,x(t),u(t),\psi_0,\psi(t)\right)
= \max_{u \in \Omega} H\left(t,x(t),u,\psi_0,\psi(t)\right) \,
\text{.}
\end{equation}
\end{description}
The Pontryagin extremal is called normal if $\psi_0 \ne 0$ and abnormal
otherwise. The celebrated \emph{Pontryagin Maximum Principle} asserts that
if $\left(x(\cdot),\,u(\cdot)\right)$ is a minimizer of the problem,
then there exists a nonzero pair $\left(\psi_0,\,\psi(\cdot)\right)$ such that
$\left(x(\cdot),\,u(\cdot),\,\psi_0,\,\psi(\cdot)\right)$ is
a Pontryagin extremal. Furthermore, the Pontryagin Hamiltonian along
the extremal is an absolutely continuous function of $t$,
\begin{equation*}
t \mapsto H\left(t,x(t),u(t),\psi_0,\psi(t)\right)
\in W_{1,1}\left([a,b];\mathbb{R}\right) \, ,
\end{equation*}
and satisfies the equality
\begin{equation}
\label{eq:dHdt}
\frac{\mathrm{d}H}{\mathrm{d}t}\left(t,x(t),u(t),\psi_0,\psi(t)\right)
= \frac{\partial H}{\partial t}\left(t,x(t),u(t),\psi_0,\psi(t)\right) \, ,
\end{equation}
for almost all $t \in [a,b]$, where on the left-hand side we have the
total derivative with respect to $t$ and on the right-hand side the partial
derivative of the Pontryagin Hamiltonian with respect to $t$
(\textrm{cf.} \cite{MR29:3316b}. See \cite{torresCM01I14E} for some
generalizations of this fact).


\section{Main Result}

To formulate a second Noether theorem in the optimal
control setting, first we need to have appropriate
notions of \emph{invariance} and \emph{Noether current}.
We propose the following ones.

\begin{definition}
A function  $C\left(t,\,x,\,u,\,\psi_0,\,\psi\right)$
which is constant along every Pontryagin extremal
$\left(x(\cdot),\,u(\cdot),\,\psi_0,\,\psi(\cdot)\right)$
of the problem,
\begin{equation}
\label{eq:def:CL}
C\left(t,\,x(t),\,u(t),\,\psi_0,\,\psi(t)\right) = k \, ,
\quad t \in [a,b] \, ,
\end{equation}
for some constant $k$, will be called a \emph{Noether current}.
The equation \eqref{eq:def:CL} is the \emph{conservation law} corresponding
to the Noether current $C$.
\end{definition}

\begin{definition}
\label{d:semi-invariant}
Let $C^m \ni p : [a,b] \rightarrow \mathbb{R}^k$ be
an arbitrary function of the independent variable. Using the notation
\begin{equation*}
\alpha(t) \doteq \left(t,x(t),u(t),p(t),\dot{p}(t),\ldots,p^{(m)}(t)\right) \, ,
\end{equation*}
we say that the optimal control problem is \emph{semi-invariant}
if there exists a $C^1$ transformation group
\begin{gather}
g : [a,b] \times \mathbb{R}^n \times \Omega \times \mathbb{R}^{k*(m+1)}
       \rightarrow
       \mathbb{R} \times \mathbb{R}^n \times \mathbb{R}^r \, , \notag \\
g\left(\alpha(t)\right) = \left(T\left(\alpha(t)\right),
X\left(\alpha(t)\right), U\left(\alpha(t)\right)\right) \, , \label{groupT}
\end{gather}
which for $p(t) = \dot{p}(t) = \cdots = p^{(m)}(t) = 0$ corresponds to the identity
transformation,
$g(t,x,u,0,0,\ldots,0) = (t,x,u)$  for all
$(t,x,u) \in [a,b] \times \mathbb{R}^n \times \Omega$,
satisfying the equations
\begin{multline}
\label{eq:invi}
\left(\lambda^0 \cdot p(t) + \lambda^1 \cdot \dot{p}(t) + \cdots
+ \lambda^m \cdot p^{(m)}(t)\right)
\frac{\mathrm{d}}{\mathrm{dt}} L\left(t,x(t),u(t)\right) \\
+ L\left(t,x(t),u(t)\right)
+ \frac{\mathrm{d}}{\mathrm{dt}} F\left(\alpha(t)\right)
= L\left(g\left(\alpha(t)\right)\right) \frac{\mathrm{d}}{\mathrm{d}t}
T\left(\alpha(t)\right) \, ,
\end{multline}
\begin{equation}
\label{eq:invii}
\frac{\mathrm{d}}{\mathrm{d}t} X\left(\alpha(t)\right)
= \varphi\left(g\left(\alpha(t)\right)\right)
\frac{\mathrm{d}}{\mathrm{d}t} T\left(\alpha(t)\right)  \, ,
\end{equation}
for some function $F$ of class $C^1$ and for some
$\lambda^0,\ldots,\lambda^m \in \mathbb{R}^k$. In this
case the group of transformations $g$ will be called a
\emph{gauge symmetry} of the optimal control problem.
\end{definition}

\begin{remark}
We use the term ``gauge symmetry'' to emphasize the fact that the
group of transformations $g$ depend on arbitrary functions. The
terminology takes origin from gauge invariance in electromagnetic
theory and in Yang-Mills theories, but it refers here to a wider
class of symmetries.
\end{remark}

\begin{remark}
The identity transformation is a gauge symmetry for any given
optimal control problem.
\end{remark}

\begin{theorem}[Second Noether theorem for Optimal Control]
\label{r:mainresult:2nt}
If the optimal control problem is semi-invariant
under a gauge symmetry \eqref{groupT}, then there exist
$k \left(m + 1\right)$ Noether currents of the form
\begin{multline*}
\psi_0 \left(\left.\frac{\partial F\left(\alpha(t)\right)}{\partial p_j^{(i)}}
\right|_{0} + \lambda_j^i L\left(t,x(t),u(t)\right)\right) + \psi(t) \cdot
\left.\frac{\partial X\left(\alpha(t)\right)}{\partial p_j^{(i)}}\right|_{0}\\
- H(t,x(t),u(t),\psi_0,\psi(t))
\left.\frac{\partial T\left(\alpha(t)\right)}{\partial p_j^{(i)}}\right|_{0}
\end{multline*}
($i=0,\ldots,m$, $j=1,\ldots,k$),
where $H$ is the corresponding Pontryagin Hamiltonian \eqref{eq:Hamiltonian}.
\end{theorem}

\begin{remark}
We are using the standard convention that $p^{(0)}(t) = p(t)$,
and the following notation for the evaluation of a term:
\begin{equation*}
\left.\left(\ast\right)\right|_{0} \doteq
\left.\left(\ast\right)\right|_{p(t) = \dot{p}(t) = \cdots = p^{(m)}(t) = 0}  \, .
\end{equation*}
\end{remark}

\begin{remark}
For the basic problem of the calculus of variations,
\textrm{i.e.}, when $\varphi = u$,
Theorem~\ref{r:mainresult:2nt} coincides with the classical
formulation of the second Noether theorem if one puts
$\lambda^i = 0$, $i=0,\ldots,m$, and $F \equiv 0$
in the Definition~\ref{d:semi-invariant},
and the transformation group $g$ is not allowed to depend on the
derivatives of the state variables (on the control variables).
In \S \ref{sec:Example} we provide an example of the calculus of
variations for which our result is applicable while previous results are not.
\end{remark}

\begin{proof}
Let $i \in \left\{0,\ldots,m\right\}$,
$j \in \left\{1,\ldots,k\right\}$,
and $\left(x(\cdot),u(\cdot),\psi_0,\psi(\cdot)\right)$ be an
arbitrary Pontryagin extremal of the optimal control problem.
Since it is assumed that to the values $p(t) = \dot{p}(t) = \cdots = p^{(m)}(t) = 0$
it corresponds the identity gauge transformation,
differentiating \eqref{eq:invi} and \eqref{eq:invii}
with respect to $p_j^{(i)}$ and then setting
$p(t) = \dot{p}(t) = \cdots = p^{(m)}(t) = 0$ one gets:
\begin{multline}
\label{eq:ds0i}
\lambda_j^i \frac{\mathrm{d}}{\mathrm{d}t} L + \frac{\mathrm{d}}{\mathrm{d}t}
\left.\frac{\partial F\left(\alpha(t)\right)}{\partial p_j^{(i)}}\right|_{0}
= \frac{\partial L}{\partial t}
      \left.\frac{\partial T\left(\alpha(t)\right)}{\partial p_j^{(i)}}\right|_{0}
      + \frac{\partial L}{\partial x} \cdot
      \left.\frac{\partial X\left(\alpha(t)\right)}{\partial p_j^{(i)}}\right|_{0} \\
      + \frac{\partial L}{\partial u} \cdot
      \left.\frac{\partial U\left(\alpha(t)\right)}{\partial p_j^{(i)}}\right|_{0}
      + L \frac{\mathrm{d}}{\mathrm{d}t}
      \left.\frac{\partial T\left(\alpha(t)\right)}{\partial p_j^{(i)}}\right|_{0} \, ,
\end{multline}
\begin{multline}
\label{eq:ds0ii}
\frac{\mathrm{d}}{\mathrm{d}t} \left.\frac{\partial X\left(\alpha(t)\right)}{\partial
p_j^{(i)}}\right|_{0}
= \frac{\partial \varphi}{\partial t}
      \left.\frac{\partial T\left(\alpha(t)\right)}{\partial p_j^{(i)}}\right|_{0}
      + \frac{\partial \varphi}{\partial x} \cdot
      \left.\frac{\partial X\left(\alpha(t)\right)}{\partial p_j^{(i)}}\right|_{0} \\
      + \frac{\partial \varphi}{\partial u} \cdot
      \left.\frac{\partial U\left(\alpha(t)\right)}{\partial p_j^{(i)}}\right|_{0}
      + \varphi \frac{\mathrm{d}}{\mathrm{d}t}
      \left.\frac{\partial T\left(\alpha(t)\right)}{\partial p_j^{(i)}}\right|_{0} \, ,
\end{multline}
with $L$ and $\varphi$, and its partial derivatives, evaluated at
$\left(t,x(t),u(t)\right)$.
Multiplying \eqref{eq:ds0i} by $\psi_0$ and \eqref{eq:ds0ii} by $\psi(t)$,
we can write:
\begin{multline}
\label{eq:joined}
\begin{split}
\psi_0 \left(\frac{\partial L}{\partial t}
      \left.\frac{\partial T\left(\alpha(t)\right)}{\partial p_j^{(i)}}\right|_{0}
      + \frac{\partial L}{\partial x} \cdot
      \left.\frac{\partial X\left(\alpha(t)\right)}{\partial p_j^{(i)}}\right|_{0}
      + \frac{\partial L}{\partial u} \cdot
      \left.\frac{\partial U\left(\alpha(t)\right)}{\partial p_j^{(i)}}\right|_{0} \right.\\
      \left. + L \frac{\mathrm{d}}{\mathrm{d}t}
      \left.\frac{\partial T\left(\alpha(t)\right)}{\partial p_j^{(i)}}\right|_{0}
- \frac{\mathrm{d}}{\mathrm{d}t}
\left.\frac{\partial F\left(\alpha(t)\right)}{\partial p_j^{(i)}}\right|_{0}
- \lambda_j^i \frac{\mathrm{d}}{\mathrm{d}t} L\right)
\end{split} \\
\begin{split}
+ \psi(t) \cdot \left(\frac{\partial \varphi}{\partial t}
      \left.\frac{\partial T\left(\alpha(t)\right)}{\partial p_j^{(i)}}\right|_{0}
      + \frac{\partial \varphi}{\partial x} \cdot
      \left.\frac{\partial X\left(\alpha(t)\right)}{\partial p_j^{(i)}}\right|_{0}
      + \frac{\partial \varphi}{\partial u} \cdot
      \left.\frac{\partial U\left(\alpha(t)\right)}{\partial p_j^{(i)}}\right|_{0}
      \right. \\
      \left. + \varphi \frac{\mathrm{d}}{\mathrm{d}t}
      \left.\frac{\partial T\left(\alpha(t)\right)}{\partial p_j^{(i)}}\right|_{0} -
      \frac{\mathrm{d}}{\mathrm{d}t}
      \left.\frac{\partial X\left(\alpha(t)\right)}{\partial p_j^{(i)}}\right|_{0}\right)
      = 0 \, .
\end{split}
\end{multline}
According to the maximality condition \eqref{eq:cm}, the function
\begin{equation*}
\psi_0 L\left(t,x(t),U\left(\alpha(t)\right)\right)
+ \psi(t) \cdot \varphi\left(t,x(t),U\left(\alpha(t)\right)\right)
\end{equation*}
attains an extremum for $p(t) = \dot{p}(t) = \cdots = p^{(m)}(t) = 0$. Therefore
\begin{equation*}
\psi_0 \frac{\partial L}{\partial u} \cdot
\left.\frac{\partial U\left(\alpha(t)\right)}{\partial p_j^{(i)}}\right|_{0}
+ \psi(t) \cdot \frac{\partial \varphi}{\partial u} \cdot
\left.\frac{\partial U\left(\alpha(t)\right)}{\partial p_j^{(i)}}\right|_{0} = 0
\end{equation*}
and \eqref{eq:joined} simplifies to
\begin{multline*}
\begin{split}
\psi_0 \left(\frac{\partial L}{\partial t}
      \left.\frac{\partial T\left(\alpha(t)\right)}{\partial p_j^{(i)}}\right|_{0}
      + \frac{\partial L}{\partial x} \cdot
      \left.\frac{\partial X\left(\alpha(t)\right)}{\partial p_j^{(i)}}\right|_{0}
      + L \frac{\mathrm{d}}{\mathrm{d}t}
      \left.\frac{\partial T\left(\alpha(t)\right)}{\partial p_j^{(i)}}\right|_{0}
      \right. \\
\left. - \frac{\mathrm{d}}{\mathrm{d}t}
\left.\frac{\partial F\left(\alpha(t)\right)}{\partial p_j^{(i)}}\right|_{0}
- \lambda_j^i \frac{\mathrm{d}}{\mathrm{d}t} L\right)
\end{split} \\
\begin{split}
+ \psi(t) \cdot \left(\frac{\partial \varphi}{\partial t}
      \left.\frac{\partial T\left(\alpha(t)\right)}{\partial p_j^{(i)}}\right|_{0}
      + \frac{\partial \varphi}{\partial x} \cdot
      \left.\frac{\partial X\left(\alpha(t)\right)}{\partial p_j^{(i)}}\right|_{0}
      + \varphi \frac{\mathrm{d}}{\mathrm{d}t}
      \left.\frac{\partial T\left(\alpha(t)\right)}{\partial p_j^{(i)}}\right|_{0}
      \right. \\
      \left. - \frac{\mathrm{d}}{\mathrm{d}t}
      \left.\frac{\partial X\left(\alpha(t)\right)}{\partial p_j^{(i)}}\right|_{0}\right)
      = 0 \, .
\end{split}
\end{multline*}
Using the adjoint system \eqref{eq:sh}
and the property \eqref{eq:dHdt}, one easily concludes
that the above equality is equivalent to
\begin{equation*}
\frac{\mathrm{d}}{\mathrm{d}t} \left(\psi_0 \left.\frac{\partial F\left(\alpha(t)\right)}{\partial p_j^{(i)}}\right|_{0} +
\psi_0 \lambda_j^i L +
\psi(t) \cdot
\left.\frac{\partial X\left(\alpha(t)\right)}{\partial p_j^{(i)}}\right|_{0}
- H \left.\frac{\partial T\left(\alpha(t)\right)}{\partial p_j^{(i)}}\right|_{0}\right)
= 0 \, .
\end{equation*}
\emph{Quod erat demonstrandum}.
\end{proof}


\section{Example}
\label{sec:Example}

Consider the following simple time-optimal problem with
$n = r = 1$ and $\Omega = (-1,1)$.
Given two points $\alpha$ and $\beta$ in the state space $\mathbb{R}$,
we are to choose an admissible pair $\left(x(\cdot),u(\cdot)\right)$,
solution of the the control equation
\begin{equation*}
\dot{x}(t) = u(t) \, ,
\end{equation*}
and satisfying the boundary conditions $x(0) = \alpha$, $x(T) = \beta$,
in such a way that the time of transfer from $\alpha$ to $\beta$ is
minimal:
\begin{equation*}
T \rightarrow \min \, .
\end{equation*}
In this case the Lagrangian is given by $L \equiv 1$ while $\varphi = u$.
It is easy to conclude that the problem is invariant under the gauge symmetry
\begin{multline*}
g\left(t,x(t),u(t),p(t),\dot{p}(t),\ddot{p}(t)\right) \\
= \left(p(t)+t,(\dot{p}(t) + 1)^2 x(t),2\ddot{p}(t) x(t) + (\dot{p}(t)+1) u(t)\right) \, ,
\end{multline*}
\textrm{i.e.}, under
\begin{equation*}
T = p(t)+t \, , \quad X = (\dot{p}(t) + 1)^2 x(t) \, , \quad U = 2\ddot{p}(t) x(t)
+ (\dot{p}(t)+1) u(t) \, ,
\end{equation*}
where $p(\cdot)$ is an arbitrary function of class $C^2\left([0,T];\mathbb{R}\right)$.
For that we choose $F = p(t)$, $\lambda^0 = \lambda^1 = \lambda^2 = 0$,
and conditions \eqref{eq:invi} and \eqref{eq:invii} follows:
\begin{align*}
L\left(T,X,U\right) \frac{\mathrm{d}}{\mathrm{d}t} T
&= \frac{\mathrm{d}}{\mathrm{d}t} \left(p(t) + t\right) = \frac{\mathrm{d}}{\mathrm{d}t} F + L(t,x(t),u(t)) \, , \\
\begin{split}
\varphi\left(T,X,U\right) \frac{\mathrm{d}}{\mathrm{d}t} T
&= \left[2\ddot{p}(t) x(t) + \left(\dot{p}(t)+1\right) u(t)\right]
\left(\dot{p}(t) + 1\right) \\
&= \frac{\mathrm{d}}{\mathrm{d}t} \left[\left(\dot{p}(t) + 1\right)^2 x(t)\right]
= \frac{\mathrm{d}}{\mathrm{d}t} X \, .
\end{split}
\end{align*}
From Theorem~\ref{r:mainresult:2nt} the two non-trivial Noether currents
\begin{gather}
\psi_0 - H \label{e:NC1} \, ,\\
2 \psi(t) x(t) \, , \label{e:NC2}
\end{gather}
are obtained.
As far as $\psi_0$ is a constant, the Noether current \eqref{e:NC1}
is just saying that the corresponding Hamiltonian $H$ is constant
along the Pontryagin extremals of the problem. This is indeed the
case, since the problem under consideration is autonomous
(\textrm{cf.} equality \eqref{eq:dHdt}). The Noether current
\eqref{e:NC2} can be understood having in mind the maximality
condition \eqref{eq:cm}
($\frac{\partial H}{\partial u} = 0 \Leftrightarrow \psi(t) = 0$).


\section{Concluding Remarks}

In this paper we provide an extension of the second Noether's theorem
to the optimal control framework. The result seems to be new
even for the problems of the calculus of variations.

Theorem~\ref{r:mainresult:2nt} admits several extensions.
It was derived, as in the original work
by Noether \cite{JFM46.0770.01,MR53:10538},
for state variables in an $n$-dimensional Euclidean space.
It can be formulated, however, in contexts where the
geometry is not Euclidean (these extensions
can be found, in the classical context, \textrm{e.g.} in
\cite{MR37:799,MR39:7485,MR43:2591}).
It admits also a
generalization for optimal control problems which are invariant
in a mixed sense, \textrm{i.e.}, which are invariant under a
group of transformations depending upon $\rho$ parameters and
upon $k$ arbitrary functions and their derivatives up to some
given order. Other possibility is to obtain
a more general version of the second Noether theorem for optimal
control problems which does not admit exact symmetries. For example,
under an invariance notion up to first-order terms in the functions
$p(\cdot)$ and its derivatives (\textrm{cf.} the quasi-invariance notion
introduced by the author in \cite{torresCM02I01E} for the first Noether theorem).
These and other questions, such as the
generalization of the first and
second Noether type theorems to constrained optimal control problems,
are under study and will be addressed elsewhere.



\end{document}